\documentclass{amsart}
\usepackage{amssymb}
\usepackage{amsfonts}
\usepackage{amsmath}
\usepackage{graphicx}
\newtheorem{thm}{Theorem}[section]
\newtheorem{corollary}[thm]{Corollary}
\newtheorem{lemma}[thm]{Lemma}
\newtheorem{proposition}[thm]{Proposition}

\theoremstyle{definition}

\newtheorem{example}[thm]{Example}

\theoremstyle{remark}

\numberwithin{equation}{section}
\def\1{{\rm (1)}}
\def\2{{\rm (2)}}
\def\3{{\rm (3)}}
\def\4{{\rm (4)}}
\def\5{{\rm (5)}}
\begin{document}

%%%%%%%%%%%%%%%%%%%%%%%%%%%%%%%%%%%%%%%%%%%%%%%%%%%%%%%%%
%%%%%%%%%%%%%%%%%%%%%%%%%%%%%%%%%%%%%%%%%%%%%%%%%%%%%%%%%
\title[On the Set of $t$-Linked Overrings of an integral
domain]{On the Set of $t$-Linked Overrings of an integral domain}

%%%%%%%%%%%%%%%%%%%%%%%%%%%%%%%%%%%%%%%%%%%%%%%%%%%%%%%%%
%%%%%%%%%%%%%%%%%%%%%%%%%%%%%%%%%%%%%%%%%%%%%%%%%%%%%%%%%
\author{A. Mimouni}

\address{Department of Mathematical Sciences,
King Fahd University of Petroleum \& Minerals, P. O. Box 5046,
Dhahran 31261, KSA}

\email{amimouni@kfupm.edu.sa}

\thanks{This work was funded by KFUPM under Project \# FT/18-2005.}

\date{}

%%%%%%%%%%%%%%%%%%%%%%%%%%%%%%%%%%%%%%%%%%%%%%%%%%%%%%%%%
%%%%%%%%%%%%%%%%%%%%%%%%%%%%%%%%%%%%%%%%%%%%%%%%%%%%%%%%%


%%%%%%%%%%%%%%%%%%%%%%%%%%%%%%%%%%%%%%%%%%%%%%%%%%%%%%%%%%%%%%%%%
%%%%%%%%%%%%%%%%%%%%%%%%%%%%%%%%%%%%%%%%%%%%%%%%%%%%%%%%%%%%%%%%%
\begin{abstract}
Let $R$ be an integral domain with quotient field $L$. An overring
$T$ of $R$ is $t$-linked over $R$ if $I^{-1}=R$ implies that
$(T:IT)=T$ for each finitely generated ideal $I$ of $R$. Let
$O_{t}(R)$ denotes the set of all $t$-linked overrings of $R$ and
$O(R)$ the set of all overrings of $R$. The purpose of this paper is
to study some finiteness conditions on the set $O_{t}(R)$.
Particularly, we prove that if $O_{t}(R)$ is finite, then so is
$O(R)$ and $O_{t}(R)=O(R)$, and if each chain of $t$-linked
overrings of $R$ is finite, then each chain of overrings of $R$ is
finite. This yields that the $t$-linked approach is more efficient
than the Gilmer's treatment in \cite{G1}. We also examine the
finiteness conditions in some Noetherian-like settings such as Mori
domain, quasicoherent Mori domain, Krull domain etc. We establish a
connection between $O_{t}(R)$ and the set of all strongly divisorial
ideals of $R$ and we conclude by a characterization of domains $R$
that are $t$-linked under all their overrings.
\end{abstract}

\maketitle

%%%%%%%%%%%%%%%%%%%%%%%%%%%%%%%%%%%%%%%%%%%%%%%%%%%%%%%%%%%%%%%%%%%%%%%%%%%%%%%%%%%%%%%%%%%%%%%%%%%%%%%%%%%%%%%%%%%%%
%%%%%%%%%%%%%%%%%%%%%%%%%%%%%%%%%%%%%%%%%%%%%%%%%%%%%%%%%%%%%%%%%%%%%%%%%%%%%%%%%%%%%%%%%%%%%%%%%%%%%%%%%%%%%%%%%%%%%
%%%%%%%%%%%%%%%%%%%%%%%%%%%%%%%%%%%%%%%%%%%%%%%%%%%%%%%%%%%%%%%%%%%%%%%%%%%%%%%%%%%%%%%%%%%%%%%%%%%%%%%%%%%%%%%%%%%%%

\begin{section}{Introduction}\label{Int}

Throughout $R$ is an integral domain (which is not a field) with
quotient field $L$. By an overring of $R$ we mean an integral domain
$T$ such that $R\subseteq T\subseteq L$. In their study of the
residually algebraic pairs of integral domains in \cite{AJ}, Ayache
and Jaballah encountered the following two conditions on the set of
overrings of an integral domain $R$:
\begin{enumerate}
\item $R$ has only finitely many overrings.
\item Each chain of distinct overrings of $R$ is finite.
\end{enumerate}
In extending results of \cite{AJ} in \cite{Ja1} and \cite{Ja2},
Jaballah asked \cite[Question 1]{Ja2} for a characterizations for
domains with condition (i), that is, domains with finitely many
overrings. In \cite{G1}, R. Gilmer labels the above two conditions
as $(FO)$ and $(FC)$ in the following meaning:
\begin{enumerate}
\item $R$ is an $FO$-domain if $R$ has finitely many overrings;
\item $R$ is an $FC$-domain if each chain of distinct overrings of $R$ is
finite.
\end{enumerate}
He completely characterizes these domains \cite[Theorem 2.14, and
Theorem 3. 4]{G1}. Gilmer's characterizations involve the integral
closure $R'$ of $R$, the conductor $(R:R')$ and the notion of $FIP$
property (stands for ``finitely many intermediate $R$-algebras")
introduced and investigated by D. Anderson, D. Dobbs and B. Mullins
in \cite{ADM}. In a recent work, A. Jaballah gives an algorithm of
how to compute the number of overrings of an integrally closed
domain with finitely many overrings \cite{Ja3}, however the
question is still open for an arbitrary integral domain.\\

According to \cite{DHLZ}, an overring $T$ of $R$ is said to be
$t$-linked over $R$ if $I^{-1}=R$ for a finitely generated ideal $I$
of $R$ implies that $(IT)^{-1}=(T:IT)=T$. The $t$-linked concept was
used in \cite[Theorem 2.10, and Corollary 2.18]{DHLZ} to find
characterizations of certain classes of $PVMD$'s analogous to
characterizations of Pr\"ufer domains due to Davis and Richman. A
domain $R$ is said to be $t$-linkative if each overring of $R$ is
$t$-linked over $R$. The class of $t$-linkative domains was
introduced (but not named) in \cite[Theorem 2.6]{DHLZ} and named, in
a detailed study, by the same authors and M. Roitman in
\cite{DHLRZ}. They also introduced the notion of super-$t$-linkative
domain, that is, a domain $R$ such that each overring of $R$ is
$t$-linkative. Let $\mathcal{O}_{t}(R)$ denote the set of all
$t$-linked overrings of $R$ and $\mathcal{O}(R)$ the set of all
overrings of $R$. Then $\mathcal{O}_{t}(R)\subseteq \mathcal{O}(R)$
and the inclusion may be strict (see \cite[Example 4.1]{DHLRZ} for a
domain $R$ such that $R'$ is not $t$-linked over $R$). Clearly an
$FO$-domain is an $FC$-domain, an $FC$-domain is super-t-linkative
(\cite[Theorem 2.14]{G1} and \cite[Corollary 2.5]{DHLRZ}) and a
super-t-linkative domain is $t$-linkative. The $t$-linkative
approach has proved its ability to be more efficient in the study of
the set of overrings of integral domains. The purpose of this paper
is to continue the investigation of some finiteness conditions on
the set of $t$-linked overrings of an integral domain. In Section 2,
we answer, in the positive, the following two questions:
\begin{enumerate}
\item If $\mathcal{O}_{t}(R)$ is finite, then is $\mathcal{O}(R)$ finite?
and do we have $\mathcal{O}_{t}(R)=\mathcal{O}(R)$?
\item If each chain of $t$-linked overrings of $R$ is finite, then is $R$ an $FC$-domain?
\end{enumerate}
As $FO$-domains and $FC$-domains have finite spectrum, we list some
results relating the (Krull) dimension of $R$ to the cardinality of
$\mathcal{O}(R)$. Section 3 is devoted to the study of
Noetherian-like settings. The main result asserts that if a Mori
domain $R$ is an $FC$-domain, then its complete integral closure
$\bar{R}$ is a Dedekind domain and the conductor $(R:\bar{R})\not
=(0)$, and so $R$ is a one-dimensional domain. Moreover, if
$A=(R:\bar{R})$ is a finitely generated ideal of $R$, then $R$ is
Noetherian (Theorem~\ref{NLS.6}).\\

Let $SD(R)$ be the set of all (nonzero) strongly divisorial ideals
of $R$ (we recall that a nonzero ideal $I$ is strongly divisorial if
$I=II^{-1}=I_{v}$) and let $\phi_{R}:SD(R)\rightarrow O_{t}(R),
I\mapsto (I:I)=I^{-1}$.  Then $\phi_{R}$ is an injective map (note
that $\phi_{R}$ is well-defined since $(I:I)$ is a $t$-linked
overring of $R$, for each strongly divisorial ideal $I$ of $R$
\cite[Proposition 2.2, (e)]{DHLZ}). In Section 4, we give necessary
and sufficient conditions for $\phi_{R}$ to be surjective and so
bijective. This leads us to compute the number of strongly
divisorial ideals for some classes of integral domains with finitely
many overrings. As an application, we prove that if $R$ is Mori and
$\phi_{R}$ is surjective, then $R$ is an $FC$-domain. We also
characterize $PVD$ (pseudo-valuation domain) $R$ for which
$\phi_{R}$ is surjective and compute its strongly divisorial ideals.
The last section deals with domains $R$ such that $R$ is $t$-linked
under all its overrings. We prove that a such domains are exactly
the one-dimensional local domains. The section closes with the study
of the transfer of this notion to the pullbacks
in order to provide original examples.\\
Throughout, we denote by $R'$ (resp. $\bar{R}$) the integral (resp.
complete integral) closure of $R$ and we use the symbol ``$\subset$"
for the strict inclusion. For a nonzero (fractional) ideal $I$ of
$R$, $I^{-1}=(R:I)=\{x\in K| xI\subseteq R\}$. The $v$- and
$t$-closures of $I$ are defined, respectively, by
$I_{v}=(I^{-1})^{-1}$ and $I_{t}=\bigcup J_{v}$, where $J$ ranges
over the set of finitely generated subideals of $I$. The ideal $I$
is said to be a $v$-ideal (or divisorial) if $I=I_{v}$, and a
$t$-ideal if $I=I_{t}$. A $t$-maximal ideal is a $t$-ideal that is
maximal for the inclusion. Finally, a local domain stands for a
domain with exactly one maximal ideal and a semilocal domain is a
domain with a finite number of maximal ideals. Unreferenced material
is standard, typically as in \cite{G2}.
\end{section}

%%%%%%%%%%%%%%%%%%%%%%%%%%%%%%%%%%%%%%%%%%%%%%%%%%%%%%%%%%%%%
\vspace{2mm}

\begin{section}{General settings}\label{GS}

We start this section by showing that if $O_{t}(R)$ is finite, then
so is $O(R)$; and $O_{t}(R)=O(R)$. We also prove that if every chain
of $t$-linked overrings is finite, then $R$ is an $FC$-domain. This
yields that the $t$-linked approach is more efficient than the
Gilmer treatment of $FO$-domains and $FC$-domains (see \cite{G1}).
As $FO$-domains and $FC$-domains have finite spectrum, so finite
dimension, we list a few results treating the relation between the
Krull dimension of an $FO$-domain $R$ and the cardinality of $O(R)$.

\begin{proposition}\label{GS.1} If $O_{t}(R)$ is finite then
$O(R)=O_{t}(R)$ and therefore $R$ is an $FO$-domain.
\end{proposition}

The proof need the following two results due to M. Zafrullah. For
the convenience of the reader, we include them with their proofs
\cite{Z1}.

\begin{proposition}\label{GS.2} Let $*$ be a star operation
of finite type on $R$. If $R$ has only a finite number of distinct
maximal $*$-ideals $P_{1}, \dots, P_{r}$, then $P_{1}, \dots, P_{r}$
are precisely the maximal ideals of $R$.
\end{proposition}

\begin{proof} Let $*$ be a star operation of finite type. Suppose that $P_{1}, \dots,
P_{r}$ are the distinct maximal $*$-ideals of $R$. Since for each
nonzero nonunit $a\in R$ the ideal $aR$ is a $*$-ideal and so is
must be contained in at least one maximal $*$-ideal, we conclude
that $R\setminus U\subseteq P_{1}\cup P_{2}\cup \dots\cup P_{r}$,
where $U$ is the group of units of $R$. Now, as each maximal ideal
$M$ consists of nonunits, we have $M\subseteq P_{1}\cup P_{2}\cup
\dots\cup P_{r}$. By the prime avoidance lemma we have $M\subseteq
P_{i}$ for some $i$. But since $M$ is maximal, we have $M=P_{i}$.
Next, as each $P_{j}$ is contained in a maximal ideal which in turn
is contained in some $P_{k}$ and since $P_{i}$ are incomparable, we
conclude that each of $P_{i}$ is a maximal ideal.
\end{proof}

\begin{corollary}\label{GS.3} If $R$ has a finite number of maximal
$t$-ideals then every maximal ideal is a $t$-ideal.
\end{corollary}

\noindent{\bf Proof of the Proposition~\ref{GS.1}} Assume that
$O_{t}(R)$ is finite. Since for each $t$-maximal ideal $M$ of $R$,
$R_{M}$ is $t$-linked over $R$, then $R$ has only finitely many
maximal $t$-ideals. By Corollary~\ref{GS.3}, every maximal ideal is
a $t$-ideal. Hence $O(R)=O_{t}(R)$ by \cite[Theorem 2.6]{DHLZ}.
$\square$\\

\begin{proposition}\label{GS.4} $R$ is an $FC$-domain if and only if
each chain of $t$-linked overrings of $R$ is finite.
\end{proposition}

\begin{proof} Assume that each chain of $t$-linked overrings of $R$
is finite. Then $R$ has finitely many maximal $t$-ideals. Indeed,
suppose on the contrary that $R$ has infinitely many maximal
$t$-ideals. Let $\{M_{i}\}_{i\geq 1}$ be an index ordered set of
$t$-maximal ideals of $R$. For each $n\geq 1$ set
$R_{n}=\displaystyle\cap_{i=1}^{i=n}R_{M_{i}}$. Then $R_{n}$ is a
$t$-linked overring of $R$ (as an intersection of $t$-linked
overrings, \cite[Proposition 2.2, (b)]{DHLZ}). Hence
$\{R_{n}\}_{n\geq 1}$ is an infinite chain of $t$-linked overrings
of $R$ \cite[Lemma 1.4]{G1}, which is absurd. So $R$ has a finite
number of $t$-maximal ideals. By Corollary~\ref{GS.3} every maximal
ideal is a $t$-ideal and by \cite[Proposition 2.2]{DHLZ},
$O(R)=O_{t}(R)$. It follows that $R$ is an $FC$-domain.
\end{proof}

It's easy to see that for any domain $R$, $1+dimR\leq |O(R)|\leq
|SSFc(R)|$, where $SSFc(R)$ is the set of all semistar operations of
finite character on $R$. In particular, if $SSFc(R)$ is finite, then
$R$ is an $FO$-domain. The next theorem restates the following
results \cite[Theorem 7]{MS} and \cite[Theorem 4.3 and Theorem
4.4]{M1} by subsituting $O(R)$ to $SSFc(R)$. The proofs are the same
with a minor changes. For the convenience of the reader, we include
them here.

\begin{thm}\label{GS.5} Let $R$ be an integral domain of
finite dimension.\\
$(1)$  $|O(R)| = 1 + dimR$ if and only if $R$ is a valuation domain
\cite[Theorem 7]{MS}.\\
$(2)$ $|O(R)|=2 + dimR$ if and only if $R$ is a local domain, $R'$
is a valuation domain and each proper overring of $R$ contains $R'$.
In this case $O(R)=\{R\}\bigcup O(R')$ \cite[Theorem 4.4]{M1}.\\
$(3)$ Assume that $R$ is not local. Then $3 + dimR\leq |O(R)|$, and
the equality holds if and only if $R$ is a Pr\"ufer domain with
exactly two maximal ideals and $Y$-graph spectrum \cite[Theorem
4.3]{M1}.
\end{thm}

\begin{proof} Set $n=dimR$ and let $(0)\subset P_{1}\subset
\dots\subset P_{n}$ be a chain of prime ideals of $R$ such that
$dimR=htP_{n}=n$ and let $\mathbb{E}=\{ L, R_{P_{1}}, \dots,
R_{P_{n}}\}$. Clearly $|\mathbb{E}|=n+1$. Now, assume that $|O(R)| =
1 + dimR$, then $O(R)=\mathbb{E}$, and therefore each overring of
$R$ is flat (as a localization of $R$). Hence $R$ is a Pr\"ufer
domain. Clearly $R$ is local with maximal ideal $P_{n}$ (and so
$R_{P_{n}}=R$). Otherwise, there is a maximal ideal $M$ of $R$ such
that $M\not =P_{n}$. Then $O(R)=\mathbb{E}\subset
\mathbb{E}\cup\{R_{M}\}$, which is absurd. Hence $R$ is local and therefore a valuation domain.\\
The converse is clear since the overrings of a valuation domain $R$
are exactly the localizations of $R$ at prime ideals.\\

$(2)$ We claim that $R$ is local. Indeed,  set $n=dimR$ and let
$(0)\subset P_{1}\subset \dots\subset P_{n}=M$ be a chain of prime
ideals of $R$ such that $dimR=htM=n$. If $N$ is a maximal ideal of
$R$ such that $N\not =M$, then $\{ R, R_{P_{1}}, \dots,
R_{P_{n-1}},R_{M}, R_{N}, L\}\subseteq O(R)$. So $|O(R)|\geq
3+n=3+dimR$, which is absurd. Hence $R$ is local. If $R$ is
integrally closed, then $R$ is a valuation domain \cite[Theorem
1.5]{G1} and so $|O(R)|=1+dimR$ by Proposition~\ref{GS.5}, which is
absurd. Hence $R\subset R'$ and therefore $1+dimR=1+dimR'\leq
|O(R')|\leq |O(R)|-1=1+dimR=1+dimR'$. So $1+dimR'=|O(R')|$ and
therefore $R'$ is a valuation domain. It follows that $O(R)=\{ R,
R_{P_{1}}, \dots, R_{P_{n-1}}, R', L\}$. Now, for each $i\in \{1,
\dots, n-1\}$, let $Q_{i}$ be a prime ideal of $R'$ such that
$Q_{i}\bigcap R=P_{i}$. Since $R'_{Q_{i}}\in O(R)$, then
$R'_{Q_{i}}=R_{P_{j}}$ for some $j\in \{1, \dots n-1\}$. Hence
$Q_{i}R'_{Q_{i}}=P_{j}R_{P_{j}}$ and therefore $P_{i}=P_{j}$. So
$i=j$ and therefore $R'\subseteq R'_{Q_{i}}=R_{P_{i}}$. Hence every
proper overring of $R$ contains $R'$ and therefore
$O(R)=O(R')\bigcup\{R\}$.\\
Conversely, assume that $R'$ is a valuation domain and
$O(R)=O(R')\bigcup\{R\}$. Then
$|O(R)|=1+|O(R')|=1+(1+dimR')=2+dimR'=2+dimR$, as desired.\\

$(3)$ Assume that $R$ is not local. Set $n=dimR$ and let $(0)\subset
P_{1}\subset \dots\subset P_{n}=M$ be a chain of prime ideals of $R$
such that $dimR=htM=n$. Let $N$ be maximal ideal of $R$ such that
$N\not =M$ and let $\mathbb{E}=\{ L, R_{P_{1}}, \dots, R_{P_{n-1}},
R_{M}, R_{N}, R\}$. Clearly $\mathbb{E}\subseteq O(R)$ and
$|\mathbb{E}|=n+3$. Hence $3 + dimR\leq |O(R)|$.  Now, assume that
$|O(R)| = 3 + dimR= |\mathbb{E}|$. Then $O(R)=\mathbb{E}$. Since $R$
is not a field and $R'$ is integral over $R$, then $R'$ cannot be a
localization of $R$. Hence $R=R'$. Then $R$ is an integrally closed
$FO$-domain. By \cite[Theorem 1.5]{G1}, $R$ is a Pr\"ufer domain
with finite spectrum, and clearly $Spec(R)=\{(0)\subset P_{1}\subset
\dots P_{n}=M, N\}$. It remains to prove that $P_{n-1}\subseteq N$.\\
Suppose that $N$ and $P_{n-1}$ are not comparable. Set
$T=R_{N}\bigcap R_{P_{n-1}}$. Then $T$ is a quasilocal Pr\"ufer
domain with exactly two maximal ideals $Q_{1}=NR_{N}\bigcap T$ and
$Q_{2}=P_{n-1}R_{P_{n-1}}\bigcap T$ (since $T$ is the intersection
of the two valuation domains $R_{N}$ and $R_{P_{n-1}}$ that are not
comparable, \cite[Propositions 1\& 2, page 412]{Bo}). Clearly
$Q_{1}\not =M$ (otherwise, $N=M$, which is absurd) and $Q_{2}\not
=M$ (otherwise, $P_{n-1}=M$, which absurd). Hence $R\not =T$. Since
all the other overrings of $R$ are valuation domains, then $O(R)=
\mathbb{E}\subset \mathbb{E}\bigcup \{T\}\subseteq O(R)$, a
contradiction. Hence $N$ and $P_{n-1}$ are
comparable and by maximality $P_{n-1}\subset N$, as desired.\\
Conversely, assume that $R$ is Pr\"ufer with exactly two maximal
ideals $M$ and $N$ and $Y$-graph spectrum, that is, $(0)\subset
P_{1}\subset \dots P_{n-1}\subseteq M\cap N$ and let $\mathbb{E}=\{
L, R_{P_{1}}, \dots, R_{P_{n-1}},R_{M}, R_{N}, R\}$. Clearly
$|\mathbb{E}|=3+n=3+dimR$ and since each overring of $R$ is a
subintersection (i.e., intersection of localizations of $R$ at some
primes), then $O(R)=\mathbb{E}$, as desired.
\end{proof}

%%%%%%%%%%%%%%%%%%%%%%%%%%%%%%%%%%%%%%%%%%%%%%%%%%%%%%%%%%%%%%%%%%
%%%%%%%%%%%%%%%%%%%%%%%%%%%%%%%%%%%%%%%%%%%%%%%%%%%%%%%%%%%%%%%%%%
For each positive integer $n\geq 1$, there exists an $n$-dimensional
local domain $R$ with $(3+dimR)$ overrings, as it's shown by the
following example.\\

\begin{example}\label{GS.8} Let $k$ be a field and $X_{1}, \dots,
X_{n+1}$ indeterminates over $k$. Set $R_{1}=k[[X_{1}^{2},
X_{1}^{5}]]$. Clearly $R_{1}$ is a one-dimensional Noetherian local
domain and it's easy to see that $O(R_{1})=\{R_{1}, k[[X_{1}^{2},
X_{1}^{3}]], k[[X_{1}]]=R_{1}', L_{1}=qf(R_{1})=k((X_{1}))\}$. So
$|O(R_{1})|=4=3+dimR_{1}$. \\
Let $V_{1}=L_{1}[[X_{2}]]=L_{1}+M_{1}$, where $M_{1}=X_{2}V_{1}$ and
set $R_{2}=R_{1}+M_{1}$. Since $R_{1}$ is local, then so is $R_{2}$
and $dimR_{2}=dimR_{1}+dimV_{1}=1+1=2$ \cite[Theorem 2.1]{BG}. Since
each overring of $R_{2}$ is comparable to $V_{1}$ \cite[Theorem
3.1]{BG} and $V_{1}$ is a $DVR$, then $O(R_{2})=\{T+M| T\in
O(R_{1})\}\cup\{L_{2}=qf(R_{2})=k((X_{1},
X_{2}))\}$. Therefore $|O(R_{2})|=|O(R_{1})|+1=4+1=5=3+dimR_{2}$.\\
By induction on $n$, assume that $R_{n}$ is a local domain with
$dimR_{n}=n$ and $|O(R_{n})|=3+dimR_{n}$. Let $L_{n}=qf(R_{n})$ and
set $V_{n}=L_{n}[[X_{n+1}]]=L_{n}+M_{n}$, where
$M_{n}=X_{n+1}V_{n}$. Set $R_{n+1}=R_{n}+M_{n}$. Then
$dimR_{n+1}=dimR_{n}+dimV_{n}=n+1$ and
$|O(R_{n+1})|=|O(R_{n})|+1=3+dimR_{n}+1=3+dimR_{n+1}$, as desired.
\end{example}

%%%%%%%%%%%%%%%%%%%%%%%%%%%%%%%%%%%%%%%%%%%%%%%%%%%%%%%%%%%%%%%%
\bigskip
%%%%%%%%%%%%%%%%%%%%%%%%%%%%%%%%%%%%%%%%%%%
We end this section by the following result which shows that the
sets $O_{t}(R)$ of all $t$-linkted overrings of $R$ and $O_{w}(R)$
of all $w$-overrings of a domain $R$ are the same and the notions of
$t$-liknative domain and $DW$-domain coincide. First we recall the
following definitions.
\begin{enumerate}
\item An overring $T$ of a domain $R$ is said to be
a $w$-overring of $R$ if $T_{w}=T$ \cite{FgMc}.
 Let
$\mathcal{O}_{w}(R)$ denote the set of all $w$-overrings of $R$.
\item A domain $R$ is said to be a $DW$-domain if each ideal of $R$ is a $w$-ideal
\cite{M}.
\end{enumerate}

%%%%%%%%%%%%%%%%%%%%%%%%%%%%%%%%%%%%%%%%%%%%%%%%%%
\begin{proposition}\label{GS.9} Let $R$ be an integral domain.Then

$(a)$ Each $w$-overring of $R$ is $t$-linked over $R$, that is,
$\mathcal{O}_{w}(R)= \mathcal{O}_{t}(R)$.

$(b)$ The following assertions are equivalent
\begin{enumerate}
\item $R$ is $t$-linkative;
\item $R$ is a $DW$-domain;
\end{enumerate}
\end{proposition}

\begin{proof}
The first part of this proposition was proved in \cite[Proposition
3.1]{C} and the second part appears in \cite{PT}. However, we give
here a simple proof. Suppose that $R$ is $t$-linkative. let $I$ be a
nonzero ideal of $R$ and let $x\in I_{w}$. Then there exists a f. g.
ideal $J$ of $R$ such that $J^{-1}=R$ and $xJ\subseteq I$. By
\cite[Theorem 2.6]{DHLZ}, $J=R$. Hence $x\in I$ and therefore
$I=I_{w}$, as desired. The converse follows also from \cite[Theorem
2.6]{DHLZ} since $I_{t}=R$ if and only $I_{w}=R$
\end{proof}
\end{section}

%%%%%%%%%%%%%%%%%%%%%%%%%%%%%%%%%%%%%%%%%%%%%%%%%%%%%%%%%%%%%%%%%
%%%%%%%%%%%%%%%%%%%%%%%%%%%%%%%%%%%%%%%%%%%%%%%%%%%%%%%%%%%%%%%%%
%%%%%%%%%%%%%%%%%%%%%%%%%%%%%%%%%%%%%%%%%%%%%%%%%%%%%%%%%%%%%%%%%
\vspace{2mm}

\begin{section}{Noetherian-like settings}\label{NLS}

Before starting this section, we recall the following useful
definitions: An integral domain $R$ is said to be:\\
$(1)$ Mori domain if $R$ satisfies the $acc$ condition on the
$v$-ideals, and seminormal if $x\in R$, for each $x\in K$ with
$x^{2}, x^{3}\in R$ (see \cite{Ba2} for more details about Mori seminormal domains).\\
$(2)$ A conducive domain if $(R:T)\not =(0)$ for each overring $T$
of $R$ with $T\subsetneq K$ (see \cite{DF}).\\
$(3)$ Almost Krull if $R_{M}$ is a Krull domain for each maximal
ideal $M$ of $R$, \cite{G2}.\\
$(4)$ A pseudo-valuation domain ($PVD$ for short) if there exists a
valuation overring $V$ of $R$ such that $Spec(R)=Spec(V)$, \cite{HH1, HH2}.\\
$(5)$ Quasi-Pr\"ufer, if $R'$ is Pr\"ufer, \cite[Corollary 6.5.14]{FHP}.\\

\begin{proposition}\label{NLS.1}
Let $R$ be a Krull domain. The following statements are equivalent
\item (i) $R$ is super-t-linkative;
\item (ii) $R$ is $t$-linkative;
\item (iii) $R$ is a Dedekind domain.
\end{proposition}

\begin{proof} $(i)\Longrightarrow (ii)$ Trivial.

$(ii)\Longrightarrow (iii)$ Assume that $R$ is $t$-linkative. By
\cite[Theorem 2.6]{DHLZ} each t-invertible ideal is invertible.
Since $R$ is Krull, then each nonzero ideal is $t$-invertible, and
so invertible. Hence $R$ is Dedekind.

$(iii)\Longrightarrow (i)$ Follows from \cite[Proposition
3.13]{DHLRZ}.
\end{proof}

\begin{corollary}\label{NLS.2} Let $R$ be an almost Krull domain.
Then $R$ is super-t-linkative if and only if $R$ is almost Dedekind.
\end{corollary}

\begin{proof} Let $M$ be a maximal ideal of $R$. Then $R_{M}$ is a
Krull domain. Since $R_{M}$ is super-t-linkative, then $R_{M}$ is a
Dedekind domain. Hence $R$ is almost Dedekind.\\
Conversely, if $R$ is almost Dedekind, then $R$ is a Pr\"ufer
domain. By \cite[Corollary 2.5]{DHLRZ}, $R$ is super-t-linkative.
\end{proof}

In \cite[Proposition 3.13]{DHLRZ}, it was proved that a Noetherian
domain $R$ is super-t-linkative if and only if $dim R\leq 1$. Our
next two results show that a one-dimensional Mori domain needs not
be super-t-linkative, and state conditions under which a
super-t-linkative Mori domain is of dimension one.

\begin{lemma}\label{NLS.3} Let $R$ be a Mori domain. If $R$ is
$t$-linkative, then every prime ideal is divisorial.
\end{lemma}

\begin{proof} Let $P$ be a prime ideal of $R$. If $htP=1$, then $P$
is divisorial, \cite[Theorem 3.1]{Ba2}, as desired. Assume that
$htP\geq 2$. Also by \cite[Theorem 3.1]{Ba2}, either $P$ is strongly
divisorial or $P^{-1}=R$. If $P^{-1}=R$, since $R$ is Mori, then
$P_{t}=P_{v}=R$, which contradicts \cite[Theorem 2.6]{DHLZ}. Hence
$P$ is divisorial, as desired.
\end{proof}

\begin{proposition}\label{NLS.4} Let $R$ be a Mori domain such that
either $(i)$ $(R:\bar{R})\not =(0)$, or $(ii)$ $R$ is seminormal. If
$R$ is super-t-linkative, then $dimR=1$.
\end{proposition}

\begin{proof} $(i)$ Assume that $(R:\bar{R})\not =(0)$. Then
$\bar{R}$ is a Krull domain, \cite[Corollary 18]{Ba1} or
\cite[Theorem 7.4]{Ba2}. Since $\bar{R}$ is super-t-linkative (as an
overring of $R$), by Proposition~\ref{NLS.1}, $\bar{R}$ is a
Dedekind domain. So $dim\bar{R}=1$. By \cite[Corollary 3.4]{BH},
$dimR=1$.\\

$(ii)$ Assume that $R$ is seminormal. Suppose that $dimR\geq 2$ and
let $M$ be a maximal ideal of $R$ such that $htM\geq 2$. Since
$R_{M}$ is a super-t-linkative Mori domain which is seminormal, then
without loss of generality, we may assume that $R$ is local with
maximal ideal $M$. Since $htM\geq 2$, by \cite[Lemma 2.5]{BH},
$B=(M:M)$ contains a nondivisorial prime ideal $Q$ such that $Q\cap
R=M$. Since $htM\geq 2$, then $M^{-1}=(M:M)$. So $B$ is a Mori
domain which is super-t-linkative and which contains a non
divisorial prime ideal, which is absurd by Lemma~\ref{NLS.3}. It
follows that $dimR=1$.
\end{proof}

The converse is not true as it's shown by the following example.
First we recall that the valuative dimension of $R$, denoted by
$dim_{v}R$, is given by $dim_{v}R=Max\{dimV| V$ valuation overring of $R$\}.\\

\begin{example}\label{NLS.5}
\end{example}
Let $k$ be a field and $X, Y$ and $Z$ indeterminates over $k$. Set
$R=k + Zk(X, Y)[[Z]]= k+M$, where $M=Zk(X, Y)[[Z]]$. Then $R$ is an
integrally closed Mori domain (in fact $R$ is a $PVD$ with maximal
ideal $M$ and associated valuation overring $V=k(X, Y)[[Z]]$).
$\bar{R}=V=k(X, Y)[[Z]]$, $(R:\bar{R})=M$. Since $dim_{v}R=3$, by
\cite[Proposition
3.12]{DHLRZ}, $R$ is not super-t-linkative.\\

From \cite[Theorem 2.14]{G1}, it's easy to see that a Noetherian
domain $R$ is an $FC$-domain if and only if $R'=\bar{R}$ is a
quasilocal Dedekind domain and $(R:R')\not =(0)$.\\

\begin{thm}\label{NLS.6} Let $R$ be a Mori domain. If $R$ is
an $FC$-domain, then $\bar{R}$ is a semilocal Dedekind domain and
$(R:\bar{R})\not =(0)$, in particular $dim_{v}R=1$. Moreover, if
$A=(R:\bar{R})$ is a finitely generated ideal of $R$, then $R$ is
Noetherian.
\end{thm}

\begin{proof} Since every $FC$-domain is super-t-linkative, by Proposition~\ref{NLS.4},
it suffices to show that $(R:\bar{R})\not =(0)$. By \cite[Theorem
2.14]{G1}, $R'$ is a Pr\"ufer domain and $R'$ is a finite
$R$-module. So the conductor $A=(R:R')\not =(0)$. Set
$T=(A_{v}:A_{v})=(AA^{-1})^{-1}$. Then $T$ is a Mori domain. Since
$A$ is an ideal of $R'$, then $R'\subseteq (A:A)\subseteq
(A_{v}:A_{v})=T$. Since $R'$ is Pr\"ufer, then so is $T$. Hence $T$
is a Dedekind domain (as a Pr\"ufer domain which is also Mori).
Clearly $(R:T)=(AA^{-1})_{v}\not =(0)$ and so $R$ and $T$ have the
same complete integral closure, that is, $\bar{R}=\bar{T}=T$. Hence
$\bar{R}$ is Dedekind and so $dim\bar{R}=1$. Since $\bar{R}$ is an
$FC$-domain (as an overring of $R$), then $\bar{R}$ has finite
spectrum. So $\bar{R}$ is semilocal, and clearly
$(R:\bar{R})=(R:T)=(AA^{-1})_{v}\not =(0)$, as desired. Since
$dimR=dim\bar{R}=1$, \cite[Corollary 3.4. (1)]{BH}, and $R'$ is a
Pr\"ufer domain, then $dim_{v}R=dim_{v}R'=dimR'=dimR=1$. Now, assume
that $A$ is a finitely generated ideal of $R$. Let $P$ be a nonzero
prime ideal of $R$. Since $dimR=1$, then $P=M$ is a maximal ideal
which is divisorial. If $MM^{-1}=R$, then $M$ is finitely generated
(as an invertible ideal). Assume that $MM^{-1}\subset R$. Then
$MM^{-1}=M$. Hence $M^{-1}=(M:M)\subseteq \bar{R}=A^{-1}$. So
$A=A_{v}\subseteq M_{v}=M$. By \cite[Theorem 2.14]{G1}, $R/A$ is
Artinian. Hence M/A is a finitely generated ideal of $R/A$. Set
$M/A=(\bar{x_{1}}, \dots, \bar{x_{n}})$, where $x_{i}\in M$ for each
$i=1, \dots, n$. Let $J=(x_{1}, \dots, x_{n})$. Clearly $M=J+A$, and
so $M$ is finitely generated. Hence every prime ideal of $R$ is
finitely generated and therefore $R$ is Noetherian.
\end{proof}

We recall that a domain $R$ is said to be quasi-coherent Mori domain
if every $t$-ideal is finitely generated \cite[page 85] {Kan}.

\begin{corollary}\label{NLS.7} Let $R$ be a quasi-coherent Mori
domain. If $R$ is an $FC$-domain, then $R$ is Noetherian.
\end{corollary}

\begin{proof} Since $A=(R:\bar{R})$ is a $v$-ideal of $R$, then $A$
is finitely generated. The conclusion follows from
Theorem~\ref{NLS.6}.
\end{proof}

According to \cite{FgMc}, a domain $R$ is said to be strong Mori
($SM$ for short) if $R$ satisfies the ascending chain conditions on
$w$-ideals. Noetherian domains are strong Mori and strong Mori
domains are Mori.

\begin{corollary}\label{NLS.8} Let $R$ be a strong Mori domain. If
$R$ is an $FC$-domain, then $R$ is a one-dimensional Noetherian
domain.
\end{corollary}

\begin{proof} Follows immediately from Theorem \ref{NLS.6} and
\cite[Corollary 1.10]{FgMc}.
\end{proof}

\end{section}
%%%%%%%%%%%%%%%%%%%%%%%%%%%%%%%%%%%%%%%%%%%%%%%%%%%%%%%%%%%%%%%%%%%%%%%%%%%
%%%%%%%%%%%%%%%%%%%%%%%%%%%%%%%%%%%%%%%%%%%%%%%%%%%%%%%%%%%%%%%%%%%%%%%%%%%
%%%%%%%%%%%%%%%%%%%%%%%%%%%%%%%%%%%%%%%%%%%%%%%%%%%%%%%%%%%%%%%%%%%%%%%%%%%
\vspace{2mm}

\begin{section}{$t$-Linked overrings and strongly divisorial ideals}\label{TLSD}

We recall that a nonzero ideal $I$ of a domain $R$ is said to be
strong (or a trace ideal) if $I=II^{-1}$ and strongly divisorial if
it is strong and divisorial, that is $I=I_{v}=II^{-1}$. Let $SD(R)$
denotes the set of all nonzero strongly divisorial ideals of $R$. By
\cite[Proposition 2.2]{DHLZ}, for each ideal $A$ of a domain $R$,
$(A_{v}:A_{v})$ is $t$-linked over $R$. In particular, if $I\in
SD(R)$, then $I^{-1}=(I:I)$ is $t$-linked over $R$. This yields an
injection map $\phi_{R}: SD(R)\longrightarrow
\mathcal{O}_{t}(R)\setminus\{L\}$, $I\mapsto I^{-1}$ (we note that
this map was introduced by V. Barucci in \cite{Ba1}). The following
proposition characterizes when $\phi_{R}$ is surjective.

\begin{proposition}\label{TLSD.1} Let $R$ be a domain.
\item 1) $\phi_{R}$ is surjective if and only if each $t$-linked
overring of $R$ is a fractional $v$-ideal of $R$.
\item 2) If $\phi_{R}$ is surjective, then:
\item i) $R$ is a conducive domain;
\item ii) Every non maximal prime ideal is strong.
\item iii) If $R$ is not local, then every maximal ideal of $R$ is not divisorial, that is,
for each maximal ideal $M$ of $R$, $M^{-1}=R$.
\end{proposition}

\begin{proof} 1) $\Longrightarrow)$ Let $T\in O_{t}(R)\setminus\{L\}$. Then there exists $I\in SD(R)$
such that $T=I^{-1}$. Hence $T$ is a fractional $v$-ideal of $R$.

$\Longleftarrow)$ Let $T\in O_{t}(R)$ and set $I=(R:T)$. Since $T$
is a fractional $v$-ideal of $R$, then $T\subseteq (I:I)\subseteq
I^{-1}=T_{v}=T$. Hence $I\in SD(R)$ and $\phi_{R}(I)=T$, as desired.

$2)$  $i)$ Since $R$ is not a field, then $dimR\geq 1$. Let $P$ be a
$t$-prime ideal of $R$. Now, for each overring $T\subset L$ of $R$,
$T_{R\setminus P}$ is $t$-linked over $R$ \cite[Proposition
2.9]{DHLZ}. Since $\phi_{R}$ is surjective, then there is $I\in
SD(R)$ such that $I^{-1}=T_{R\setminus P}$. Hence
$I=I_{v}=(R:I^{-1})=(R:T_{R\setminus P})\subseteq (R:T)$. So
$(R:T)\not=(0)$ and therefore $R$ is conducive.

ii) Let $P$ be a nonmaximal prime ideal of $R$. Since $R_{P}\in
O_{t}(R)\setminus\{L\}$, then $R_{P}=I^{-1}$ for some $I\in SD(R)$.
Hence $IR_{P}=II^{-1}=I$. Since $P$ is not maximal, then $R\subset
R_{P}=I^{-1}$. Hence $I\subseteq P$. So $P^{-1}\subseteq
I^{-1}=R_{P}$ and then $PP^{-1}\subseteq PR_{P}$. Since
$PP^{-1}\subseteq R$, then $PP^{-1}\subseteq PR_{P}\cap R=P$. Hence
$PP^{-1}=P$, as desired.\\

iii) Assume that $R$ is not local. Let $M$ be a maximal ideal of $R$
and $I\in SD(R)$ such that $I^{-1}=R_{M}$. Then $I=II^{-1}=IR_{M}$.
Since $R\subset R_{M}$, then $I\subseteq M$. So $M^{-1}\subseteq
I^{-1}=R_{M}$. Since $M^{-1}\subseteq R_{N}$ for each maximal ideal
$N\not =M$, then $M^{-1}\subseteq R_{M}\cap
(\displaystyle\cap_{N\not =M}R_{N})=R$. Hence $M^{-1}=R$, as
desired.
\end{proof}

\begin{corollary}\label{TLSD.2} Let $R$ be a Mori domain. If $\phi_{R}$
is surjective, then $R$ is a one-dimensional conducive local domain
and either $(i)$ $R$ is a $DVR$, or $(ii)$ the maximal ideal $M$ of
$R$ is strongly divisorial and each proper overring of $R$ contains
$M^{-1}$.
\end{corollary}

\begin{proof} By Proposition~\ref{TLSD.1}, $R$ is conducive. So
$\bar{R}$ is a one-dimensional valuation domain (it easy to see that
the only overrings of $\bar{R}$ are $\bar{R}$ and $L=qf(R)$). Since
$(R:\bar{R})\not =(0)$, then  $dimR=1$ \cite[Corollary 3.4]{BH}.
Since $R$ is Mori, then $\bar{R}$ is Krull and therefore $\bar{R}$
is a $DVR$. Set $Spec(\bar{R})=\{(0)\subset N\}$. Now, if $Q$ and
$M$ are maximal ideals of $R$, since $htQ=htM=dimR=1$, by
\cite[Proposition 1.1]{BH}, there is $Q'$ and $M'$ in
$Spec(\bar{R})$ such that $Q'\cap R=Q$ and $M'\cap R=M$. Since
$Spec(\bar{R})=\{(O)\subseteq N\}$, then $Q'=M'=N$ and therefore
$Q=M$. Hence $R$ is local with maximal ideal $M=N\cap R$.\\
1) If $MM^{-1}=R$, then $R$ is a $DVR$.\\
2) If $MM^{-1}=M$, then set $T=M^{-1}=(M:M)$. Clearly $T$ is a
proper overring of $R$ (since $M$ is a $v$-ideal of $R$). Let $S$ be
a proper overring of $R$. Since $dimR=1$, then $O(R)=O_{t}(R)$. So
$S$ is $t$-linked over $R$. But, since $\phi_{R}$ is surjective,
then there is $I\in SD(R)$ such that $S=I^{-1}$. Since $S$ is a
proper overring of $R$, then $I\subset R$. So $I\subseteq M$ (since
$R$ is local with maximal ideal $M$). Hence $T=M^{-1}\subseteq
I^{-1}=S$, as desired.
\end{proof}

\begin{corollary}\label{TLSD.3} Let $R$ be a Mori domain. If $\phi_{R}$ is
surjective, then $R$ is an $FC$-domain.
\end{corollary}

\begin{proof} By Corollary~\ref{TLSD.2}, $R$ is conducive and
$dimR=1$. So $O(R)=O_{t}(R)$. Set $A=(R:\bar{R})\not =(0)$ and let
$\{R_{n}\}_{n\geq 1}$ be a chain of overrings of $R$. Since
$\phi_{R}$ is surjective, then for each $n\geq 1$, there exists
$I_{n}\in SD(R)$ such that $I_{n}^{-1}=R_{n}$. So for each $n\geq
1$, $R_{n}=I_{n}^{-1}=(I_{n}:I_{n})\subseteq \bar{R}=A^{-1}$. Hence
$A=A_{v}\subseteq (I_{n})_{v}=I_{n}$. So $A\subseteq \displaystyle
\bigcap_{n\geq 1}I_{n}$ and therefore $\displaystyle \bigcap_{n\geq
1}I_{n}\not =(0)$. Since $I_{n}^{-1}=R_{n}\subseteq
R_{n+1}=I_{n+1}^{-1}$, then $I_{n+1}=(I_{n+1})_{v}\subseteq
(I_{n})_{v}=I_{n}$. Hence $\{I_{n}\}_{n\geq 1}$ is a decreasing
chain of $v$-ideals of $R$ with nonzero intersection. Since $R$ is
Mori, then $\{I_{n}\}_{n\geq 1}$ stabilizes, and therefore
$\{R_{n}\}_{n\geq 1}$ stabilizes. It follows that $R$ is an
$FC$-domain.
\end{proof}

\begin{proposition}\label{TLSD.4} $1)$ Let $R$ be a valuation
domain. Then $\phi_{R}$ is surjective. Moreover, if $dimR=n\geq 1$
is finite, then $R$ has exactly $n$ strongly divisorial ideals.\\
$2)$ If $R$ is completely integrally closed, then $\phi_{R}$ is
surjective if and only if $R$ is a one-dimensional valuation domain.\\
$3)$ If $R$ is a Pr\"ufer domain, then $\phi_{R}$ is surjective if
and only if $R$ is a valuation domain.
\end{proposition}

\begin{proof}$1)$ Assume that $R$ is a valuation domain. Then
$O_{t}(R)=O(R)$. Let $T\in O(R)\setminus\{L\}$. If $T=R$, then
$T=\phi_{R}(R)$. Assume that $R\subset T$. Then $T=R_{Q}$ for some
nonzero nonmaximal prime ideal $Q$ of $R$. By \cite[Theorem 3.8 and
Proposition 3.10]{HuP}, $Q\in SD(R)$ and $\phi_{R}(Q)=R_{Q}=T$.
Moreover, if $dimR=n$ is finite, then $SD(R)=\{$ all nonzero
nonmaximal prime ideals of $R\}\bigcup \{R\}$.\\

$2)$ Assume that $R$ is completely integrally closed and $\phi_{R}$
is surjective. It's clear that $SD(R)=\{R\}$. Since $\phi_{R}$ is
surjective, then $O(R)=\{R, L\}$ and therefore $R$ is
one-dimensional valuation domain.\\

$3)$ Assume that $R$ is Pr\"ufer and $\phi_{R}$ is surjective.
Clearly $O(R)=O_{t}(R)$ since the $t$-operation on a Pr\"ufer domain
is trivial. It suffices to prove that $R$ is local. If not, let $M$
be a maximal ideal of $R$. Since $\phi_{R}$ is surjective, then
there exists a proper ideal $I\in SD(R)$ such that $I^{-1}=R_{M}$.
Since $IR_{M}=II^{-1}=I$, then $I\subseteq M$. Let $N$ be a maximal
ideal of $R$ such that $N\not =M$. If $I\not\subseteq N$, then
$R_{M}=I^{-1}\subseteq R_{N}$. So $N\subseteq M$ and hence $M=N$,
which is absurd. It follows that $I\subseteq N$. Therefore $I$ is
contained in all maximal ideals of $R$. Now, let $P$ be a minimal
prime over $I$. By \cite[Theorem 3.2]{HuP}, $R_{M}=I^{-1}\subseteq
R_{P}$. Hence $P\subseteq M$. Since $R$ is a Pr\"ufer domain, then
$Spec(R)$ is treed. So $Min(I)$, the set of all minimal primes over
$I$, is reduced to one element, say $Min(I)=\{P\}$. Let $N$ be a
maximal ideal of $R$. Since $I\subseteq N$ and $Min(I)=\{P\}$, then
$P\subseteq N$. We claim that $I^{-1}=R_{P}$. Indeed, by
\cite[Theorem 3.2]{HuP}, $I^{-1}\subseteq R_{P}$. Conversely, if
there exists $x\in R_{P}$ and $a\in I$ such that $xa\not\in R$, then
$xa\not\in R_{N}$ for some maximal ideal $N$ of $R$. Since $R_{N}$
is a valuation domain, then $(xa)^{-1}\in R_{N}$. Since $P\subseteq
N$, then $R_{N}\subseteq R_{P}$. So $(xa)^{-1}\in R_{P}$. Hence
$1=(xa)^{-1}(xa)\in IR_{P}\subseteq PR_{P}$, which is absurd. It
follows that $R_{M}=I^{-1}=R_{P}$ and therefore $P=M$. Hence $M$ is
the unique minimal prime over $I$ and $M$ is contained in all
maximal ideals of $R$, which is a contradiction. It follows that $R$
is local and therefore $R$ is a valuation domain.
\end{proof}

\begin{thm}\label{TLSD.4} Let $R$ be a $PVD$ domain, $V$ its associated
valuation overring, $M$ its maximal ideal, $k=R/M$ its residue
field, $K=V/M$, and suppose that $R\subset V$. The following
statements are equivalent.

\item i) $\phi_{R}$ is surjective;
\item ii) The extension $R\subset V$ is minimal, i.e. there are no rings properly between $R$ and $V$;
\item iii) $K$ is algebraic over $k$ and for each $\alpha\in K\setminus k$,
$K=k(\alpha)$;
\item iv) $O(R)=\{R\}\bigcup O(V)$.\\
Moreover, if $dimV=n$ is finite, then $R$ is an $FO$-domain and $R$
has exactly $n+1$ strongly divisorial ideals $\{P_{1}, P_{2},\dots,
P_{n-1}, M, R\}$, where $P_{1}, P_{2},\dots, P_{n-1}, M$ are the
nonzero prime ideals of $R$.
\end{thm}

\begin{proof}  We first note that $O_{t}(R)=O(R)$ since $R$ is local with
maximal ideal $M$ and $M^{-1}=V$. \\

$i)\Longrightarrow ii)$ Assume that $\phi_{R}$ is surjective. Let
$T$ be an overring of $R$ such that $R\subset T\subseteq V$. By $i)$
there exists $I\in SD(R)$ such that $I^{-1}=T$. Since $R\subset T$,
then $I\subset R$ and since $R$ is local, then $I\subseteq M$. So
$V=M^{-1}\subseteq I^{-1}=T\subseteq V$. Hence $V=T$ and therefore
the extension
$R\subset V$ is minimal.\\

$ii)\Longrightarrow iii)$ Let $\alpha\in K\setminus k$ and set
$T=\psi^{-1}(k[\alpha])$. Since $R\subset T\subseteq V$, by $ii)$,
$T=V$. So $k[\alpha]=K$, as desired.\\

$iii)\Longrightarrow iv)$ It's well-known that each overring of $R$
is comparable to $V$. So $O(R)=[R, V[\bigcup O(V)$, where $[R, V[$
is the set of all overrings $T$ of $R$ such that $R\subseteq
T\subset V$. Hence $O(R)=\{R\}\bigcup O(V)$ since $[R, V[=\{R\}$.\\

$iv)\Longrightarrow i)$ Assume that $O(R)=\{R\}\bigcup O(V)$.
Clearly $R=\phi_{R}(R)$, $V=M^{-1}=\phi_{R}(M)$ and $M\in SD(R)$.
Let $T$ be an overring of $R$. By $iv)$ and without loss of
generality, we may assume that $T=V_{P}$, where $P$ is a nonmaximal
prime ideal of $V$. Since $R$ and $V$ have the same spectrum, $P$ is
a prime ideal of $R$ and $V_{P}=(P:P)=P^{-1}$. So $P_{v}\in SD(R)$
and $T=V_{P}=(P_{v}:P_{v})=\phi_{R}(P_{v})$. It follows that
$\phi_{R}$ is surjective.
\end{proof}
%%%%%%%%%%%%%%%%%%%%%%%%%%%%%%%%%%%%%%%%%%%%%%%%%%%%%%%%%%%%%%%%%%%%%%
%%%%%%%%%%%%%%%%%%%%%%%%%%%%%%%%%%%%%%%%%%%%%%%%%%%%%%%%%%%%%%%%%%%%%%
%%%%%%%%%%%%%%%%%%%%%%%%%%%%%%%%%%%%%%%%%%%%%%%%%%%%%%%%%%%%%%%%%%%%%%

It's clear that any $(n-1)$-dimensional valuation domain, $n\geq 2$,
has exactly $n$ overrings. Also it's easy to see that a domain $R$
has exactly two overrings if and only if $R$ is one-dimensional
valuation domain. The above theorem leads us to construct a
$PVD$ domain $R_{n}$ which is not a valuation domain, and which has exactly $n$ overrings, for any positive integer $n\geq 3$.\\

\begin{example}\label{TLSD.5}
Let $n$ be a positive integer with $n\geq 3$, $\mathbb{Q}$ be the
field of rational numbers and $X_{1}, \dots, X_{n-2}$ indeterminates
over $\mathbb{Q}$. Our aim is to construct a descending chain of
valuation domains $V_{1}\supsetneqq V_{2}\supsetneqq \dots
\supsetneqq V_{n-2}$ with maximal ideals $M_{1}\subset M_{2}\subset
\dots\subset M_{n-2}$ (respectively) such that $dimV_{i}=i$ for each
$i\in \{1, \dots, n-2\}$ and $V_{n-2}=\mathbb{Q}(\sqrt{2})+
M_{n-2}$.\\
$\centerdot$ For $n=3$, just set
$V_{1}=\mathbb{Q}(\sqrt{2})[[X_{1}]]=\mathbb{Q}(\sqrt{2})+M_{1}$,
where $M_{1}=X_{1}V_{1}$.\\
$\centerdot$ For $n=4$, set
$V_{1}=\mathbb{Q}(\sqrt{2})((X_{1}))[[X_{2}]]=\mathbb{Q}(\sqrt{2})((X_{1}))+M_{1}$,
where $M_{1}=X_{2}V_{1}$, and
$V_{2}=\mathbb{Q}(\sqrt{2})[[X_{1}]]+M_{1}=\mathbb{Q}(\sqrt{2})+M_{2}$,
where $M_{2}=X_{1}V_{2}$.\\
$\centerdot$ For $n=5$, set $V_{1}=\mathbb{Q}(\sqrt{2})((X_{1},
X_{2}))[[X_{3}]]=\mathbb{Q}(\sqrt{2})((X_{1}, X_{2}))+M_{1}$, where
$M_{1}=X_{3}V_{1}$,
$V_{2}=\mathbb{Q}(\sqrt{2})((X_{1}))[[X_{2}]]+M_{1}=\mathbb{Q}(\sqrt{2})((X_{1}))+M_{2}$,
where $M_{2}=X_{2}V_{2}$, and
$V_{3}=\mathbb{Q}(\sqrt{2})[[X_{1}]]+M_{2}=\mathbb{Q}(\sqrt{2})+M_{3}$,
where $M_{3}=X_{1}V_{3}$.\\
Iterating this process, we construct the desired chain as follows:\\
Let $V_{1} = \mathbb{Q}(\sqrt{2})((X_{1}, \dots,
X_{n-3}))[[X_{n-2}]]=\mathbb{Q}(\sqrt{2})((X_{1}, \dots, X_{n-3})) +
M_{1}$, where $M_{1}=X_{n-2}V_{1}$.\\
$V_{2} = \mathbb{Q}(\sqrt{2})((X_{1}, \dots, X_{n-4}))[[X_{n-3}]] +
M_{1}=\mathbb{Q}(\sqrt{2})((X_{1}, \dots, X_{n-4})) + M_{2}$, where
$M_{2}=X_{n-3}V_{2}$.\\
$V_{3} = \mathbb{Q}(\sqrt{2})((X_{1}, \dots, X_{n-5}))[[X_{n-4}]] +
M_{2}=\mathbb{Q}(\sqrt{2})((X_{1}, \dots, X_{n-5})) + M_{3}$, where
$M_{3}=X_{n-4}V_{3}$.\\
$V_{i}=\mathbb{Q}(\sqrt{2})((X_{1}, \dots, X_{n-i-2})) +M_{i}$,
where $M_{i}=X_{n-i-1}V_{i}$. Now, let
$R_{n} = \mathbb{Q} + M_{n-2}$. Then:\\
$|O(R_{n})|=1+ |O(V_{n-2})|= 1+ (1+dimV_{n-2})=1+1+n-2=n$, as
desired.
\end{example}

\end{section}
%%%%%%%%%%%%%%%%%%%%%%%%%%%%%%%%%%%%%%%%%%%%%%%%%%%%%%%%
\vspace{2mm}

\begin{section}{Domains that are $t$-linked under its Overrings}\label{DTUO}

According to \cite{AZ}, A domain $R$ is said to be $t$-linked under
$T$ if for each finitely generated ideal $I$ of $R$, $(T:IT)=T$
implies that $I^{-1}=R$. In some sense the notion of ``$t$-linked
under" is the opposite of the $t$-linkdness. It was introduced by
Anderson and Zafrullah in \cite{AZ} to characterize almost
B\'{e}zout domains. In this short section, we will prove some
results related to this notion. First we not that if $T=L$, then
$(L:IL)=L$ for every nonzero ideal $I$ of $R$, however $I^{-1}$ needs
not be equal to $R$. It turns that if $R$ is $t$-linked under its quotient field $L$, then $R=L$.
In this view, and in accordance with our hypothesis in the introduction we always assume that $R\subseteq T\subset L$.\\

\begin{proposition}\label{DTUO.1} Let $R$ be an integral domain and
$Q$ a nonzero prime ideal of $R$. If $R$ is $t$-linked under
$R_{Q}$, then $R$ is local with maximal ideal $Q$.
\end{proposition}

\begin{proof} Let $Q$ be a nonzero prime ideal of $R$ and let $a\in
R\setminus Q$. Set $I=aR$. Since $(R_{Q}:IR_{Q})=R_{Q}$ and $R$ is
$t$-linked under $R_{Q}$, then $a^{-1}R=I^{-1}=R$. Then $a^{-1}\in
R$ and therefore $R$ is local with maximal ideal $Q$.
\end{proof}

\begin{proposition}\label{DTUO.2} Let $R$ be an integral domain. Then $R$ is
$t$-linked under each overring $T$ of $R$ if and only if $R$ is a
one-dimensional local domain.
\end{proposition}

\begin{proof} By Proposition~\ref{DTUO.1}, $R$ is local with maximal
ideal $M$, $dimR=1$ and $Spec(R)=\{(0)\subset M\}$.\\
Conversely, assume that $R$ is local with maximal ideal $M$ and
$dimR=1$. If $|O(R)|=2$, then $R$ is the only overring of $R$
properly contained in $L$, and trivially $R$ is $t$-linked under
itself. Assume that $|O(R)|> 2$. Let $T$ be an overring of $R$ which
is not a field and let $Q$ a prime $t$-ideal of $T$ (such a prime
$t$-ideal exists since $T$ is not a field and it suffices to
consider any minimal prime over a proper principal ideal). Since $R$
is local and $dimR=1$, then $Q\cap R=M$. Now, let $I$ be a finitely
generated ideal of $R$ such that $(T:IT)=T$. If $R\subset I^{-1}$,
then $I_{v}$ is a proper ideal of $R$. So $I_{v}\subseteq M$. Hence
$IT\subseteq I_{v}T\subseteq MT\subseteq Q$. Therefore
$T=(IT)_{t_{1}}\subseteq Q_{t_{1}}=Q$ (here $t_{1}$ denotes the
$t$-operation with respect to $T$),  which is absurd since $Q$ is
$t$-ideal. Hence $I^{-1}=R$ and therefore $R$ is $t$-linked under
$T$.
\end{proof}

\begin{proposition}\label{DTUO.3} Let $R$ be a Pr\"ufer domain. If $T$ is
an overring of $R$ such that $R$ is $t$-linked under $T$, then
$R=T$.
\end{proposition}

\begin{proof} Let $T$ be an overring of $R$ such that $T\subset L$.
Assume that $R$ is $t$-linked under $T$. By \cite[Theorem 26.1]{G2}
or \cite[Theorem 1.1.2]{FHP}, $T=\bigcap \{R_{P}| P\in \Omega\}$,
where $\Omega$ is the set of all prime ideals $P$ of $R$ such that
$PT\subset T$. Let $Q$ be a nonzero prime ideal of $R$. If $QT=T$,
then $1=\displaystyle\sum_{i=1}^{i=n}a_{i}x_{i}$ where $a_{i}\in Q$
and $x_{i}\in T$ for each $i\in \{1, \dots, n\}$. Set
$I=\displaystyle\sum_{i=1}^{i=n}a_{i}R$. Then $I$ is a finitely
generated ideal of $R$ and $IT=T$. So $(T:IT)=T$. Since $R$ is
$t$-linked under $T$, then $I^{-1}=(R:I)=R$. But since $R$ is a
Pr\"ufer domain, then $I$ is invertible. So $R=II^{-1}=IR=I$ and
hence $R=I\subseteq Q$, which is a contradiction. Hence $QT\subset
T$. So $Q\in \Omega$ and therefore $\Omega = Spec(R)\setminus
\{(0)\}$. It follows that $T=\bigcap \{R_{P}| P\in \Omega\} =\bigcap
\{R_{P}| P\in Spec(R)\setminus\{(0)\}\}=R$, as desired.
\end{proof}

\begin{proposition}\label{DTUO.4} Let $R$ be an integral domain. If
$R'$ is Pr\"ufer or $dimR=1$, then $R$ is $t$-linked under $R'$.
\end{proposition}

\begin{proof} If $R=R'$, nothing to prove. So we may assume that $R\subset R'$.
Let $I$ be a nonzero finitely generated ideal of $R$ such that
$(R':IR')=R'$ and suppose that $I\subset R$. Let $P$ be a prime
ideal of $R$ such that $I\subseteq P$ and let $Q$ be a prime ideal
of $R'$ such that $Q\cap R=P$ (such a prime ideal exists since $R'$
is integral over
$R$, so the extension $R\subset R'$ is Lying Over). \\
Assume that $R'$ is Pr\"ufer. Since $IR'$ is a finitely generated
ideal, then $R'=IR'(R':IR')=IR'\subseteq PR'\subseteq Q$, which is
absurd. Hence $I=R$ and therefore $I^{-1}=R$, as desired.\\
Now, if $dimR=1$. Then $dimR'=1$. So $htQ=1$ and therefore $Q$ is a
$t$-prime ideal of $R'$. But $(R':IR')=R'$ implies that
$R'=(IR')_{v'}=(IR')_{t}\subseteq (Q)_{t'}=Q$, where $t'$ and $v'$
are the $t$- and $v$-operations with respect to $R'$. This yields a
contradiction. Hence, also in this case, $I=R$ and therefore
$(R:I)=R$, as desired.
\end{proof}

It's not the case that $R$ is always $t$-linked under $R'$ as it's
shown by the following example.

\begin{example}\label{DTUO.5} Let $\mathbb{Q}$ be the field of
rational numbers and $X$ and $Y$ indeterminates over $\mathbb{Q}$.
Set $T=\mathbb{Q}(\sqrt{2})[[X, Y]]=\mathbb{Q}(\sqrt{2})+ M$ and let
$R=\mathbb{Q}+ M$.  Clearly $R$ is a Noetherian domain which is
local with maximal ideal $M$ and $dimR=2$ \cite[Theorem 4.12]{GH}.
Since $(R:T)=M$, then $R'=\bar{R}=\bar{T}=T$. However, $R$ is not
$t$-linked under $T$ since $M$ is a finitely generated ideal of $R$,
$(T:MT)=(T:M)=T$ and $M^{-1}=T$ (or by Theorem~\ref{DTUO.6} below,
since $M$ is not a $t$-ideal of $T$).
\end{example}

%%%%%%%%%%%%%%%%%%%%%%%%%%%%%%%%%%%%%%%%%%%%%%%%%%%%%%%%%%%%%%%
\bigskip
%%%%%%%%%%%%%%%%%%%%%%%%
In the literature, a few examples satisfying this notion are known.
As pullbacks are known as a source for examples and counterexamples,
the next results characterize this notion under pullbacks.\\
Let $T$ be an integral domain which is not a field, $M$ a maximal
ideal of $T$, $K=T/M$ its residue field, $\phi:T\rightarrow K=T/M$
the canonical projection and $D$ a subring of $K$. Let $R$ be the
pullback of the following diagram

\[\begin{array}{ccl}
R            & \longrightarrow                 & D\\
\downarrow   &                                 & \downarrow\\
T            & \stackrel{\phi}\longrightarrow  & K=T/M
\end{array}\]

We assume that $R\subset T$ and we refer to this diagram as diagram
of type $(\square)$.\\

\begin{thm}\label{DTUO.6} For the diagram of type $(\square)$,
$R$ is $t$-linked under $T$ if and only $D=k$ is a field and $M$ is
a $t$-ideal of $T$.
\end{thm}

\begin{proof} Let $0\not =d\in D$,  $J=dD$ and $I=\phi^{-1}(J)$. By
\cite[Corollary 1.7]{FG}, $I$ is a finitely generated ideal of $R$.
Since $M\subset I$ and $M$ is a maximal ideal of $T$, then $IT=T$.
Hence $(T:IT)=T$. Since $R$ is $t$-linked under $T$, then
$I^{-1}=R$. By \cite[Proposition 6]{HKLM1}, $d^{-1}D=J^{-1}=D$.
Hence $d^{-1}\in D$ and therefore $D$ is a field.

Let $t_{1}$ denote the $t$-operation on $T$. Suppose that $M$ is not
a $t$-ideal of $T$. Then $M_{t_{1}}=T$. So there exists a finitely
generated subideal $A=\displaystyle\sum_{i=1}^{i=r}a_{i}T$ of $M$
such that $(T:A)=T$. Set $B=\displaystyle\sum_{i=1}^{i=r}a_{i}R$.
Then $B$ is a f.g. subideal of $M$ with $(T:BT)=(T:A)=T$. Since $R$
is $t$-linked under $T$, then $B^{-1}=R$. Hence $R=B_{v}\subseteq
M_{v}=M$ \cite[Corollary 2]{HKLM1}, which is absurd. It follows that
$M$ is a $t$-ideal of $T$.

Conversely, assume that $D=k$ is a field and $M$ is a $t$-maximal
ideal of $T$. Let $I$ be a f.g. ideal of $R$ such that $(T:IT)=T$.
Then $I^{-1}=(R:I)\subseteq (T:IT)=T$ and so $M=(R:T)\subseteq
I_{v}=I_{t}$. Since $D=k$ is a field, then $M$ is a maximal ideal of
$R$. Hence either $M=I_{v}=I_{t}$ or $I_{v}=I_{t}=R$. Suppose that
$M=I_{v}$. Then $I\subseteq M$. So $IT\subseteq M$ and therefore
$T=(IT)_{v_{1}}=(IT)_{t_{1}}\subseteq M_{t_{1}}=M$, which is absurd.
Hence $I_{v}=R$ and therefore $I^{-1}=R$. It follows that $R$ is
$t$-linked under $R$.
\end{proof}

%%%%%%%%%%%%%%%%%%%%%%%%%%%%%%%%%%%%%%%%%%%%%%%%%%%%%%%%%%%%%%%%%%
\begin{corollary}\label{DTUO.7}
$(1)$ Let $R$ be a $PVD$ (pseudo-valuation domain) and $V$ its
associated valuation overring. Then $R$ is $t$-linked under $V$.\\
$(2)$ Let $R=k+XK[X]$, where $k\subset K$ is an extension of fields.
Then $R$ is $t$-linked under $K[X]$.
\end{corollary}

%\vspace{2mm} {\bf Acknowledgment} I would like to express my sincere
%thanks to the referee for his/her many helpful suggestions and
%comments.

\end{section}

%%%%%%%%%%%%%%%%%%%%%%%%%%%%%%%%%%%%%%%%%%%%%%%%%%%%%%%%%
%%%%%%%%%%%%%%%%%%%%%%%%%%%%%%%%%%%%%%%%%%%%%%%%%%%%%%%%%
%\bibliographystyle{amsplain}
\bigskip
%%%%%%%%%%%%%%%%%%%%%%%%%%%%%%%%%%%%%%%%%%%%%%%%%%%%%%%%%
%%%%%%%%%%%%%%%%%%%%%%%%%%%%%%%%%%%%%%%%%%%%%%%%%%%%%%%%%

\begin{thebibliography}{99}
\bibitem{ADM}   D. D. Anderson, D. Dobbs and B. Mullins, \emph{The primitive element theorem for commutative algebras}, Houston J. Math. \textbf{25} (1999) 603-623.
\bibitem{AZ}    D. D. Anderson and M. Zafrullah, \emph{Almost B\'ezout domains}, J. Algebra \textbf{142} (1991) 285--309.
\bibitem{AJ}    A. Ayache and A. Jaballah, \emph{Residually algebraic pairs of rings}, Math. Zeit. \textbf{225} (1997) 49--65.
\bibitem{Ba1}   V. Barucci, \emph{Strongly divisorial ideals and complete integral closure of an integral domain}, J. Algebra \textbf{99} (1986) 132--142.
\bibitem{Ba2}   V. Barucci, \emph{Mori domains}, in ``Non-Noetherian Commutative Ring Theory" Mathematics and its Applications, Kluwer (S. Chapman and S Galz ed.) Vol.520 (2000) 57--73.
\bibitem{BH}    V. Barucci and E. Houston, \emph{On the prime spectrum of a Mori domain}, Comm. Algebra \textbf{24} (1996) 3599--3622.
\bibitem{BG}    E. Bastida and R. Gilmer, \emph{Overrings and divisorial ideals of rings of the form D+M}, Michigan Math. J. \textbf{20} (1973) 79--95.
\bibitem{Bo}    N. Bourbaki, \emph{Elements of Mathematics, Commutative Algebra}, Hermann, 1972.
\bibitem{DF}    D. Dobbs and R. Fedder, \emph{Conducive integral domains}, J. Algebra \textbf{86} (1984) 494--510.
\bibitem{DHLRZ} D. Dobbs, E. Houston, T. Lucas, M. Roitman and M. Zafrullah, \emph{On $t$-linked overrings}, Comm. Algebra \textbf{20} (1992) 1463--1488.
\bibitem{DHLZ}  D. Dobbs, E. Houston, T. Lucas and M. Zafrullah, \emph{$t$-linked overrings and Pr\"ufer $v$-multiplication domains}, Comm. Algebra \textbf{17} (1989) 2835--2852.
\bibitem{C}     G. W. Chang, \emph{Strong Mori domains and the ring $D[X]_{N_{v}}$}, J. Pure Appl. Algebra \textbf{197} (2005) 279--292.
\bibitem{FgMc}  W. Fanggui and R. L. McCasland, \emph{On strong Mori domains}, J. Pure Appl. Algebra \textbf{135} (1999) 155--165.
\bibitem{FG}    M. Fontana and S. Gabelli, \emph{On the Class group and the Local Class Group of a Pullback}, J. Algebra \textbf{181} (1996) 803--835.
\bibitem{FHP}   M. Fontana, J. Huckaba and I. Papick, \emph{Pr\"ufer domains} Marcel Dekker, \textbf{203} New York, 1997.
\bibitem{FHPR}  M. Fontana, J. Huckaba, I. Papick and M. Roitman, \emph{Pr\"ufer domains and endomorphisms of their ideals}, J. Algebra \textbf{157} (1993) 489--516.
\bibitem{GH}    S. Gabelli and E. Houston, \emph{Coherentlike Conditions in Pullbacks}, Michigan Math. J. \textbf{44} (1997) 99--122.
\bibitem{HH1}   J. Hedstrom and E. Houston, \emph{Pseudo-valuation domains}, Pacific J. Math. \textbf{75} (1978) 137--147.
\bibitem{HH2}   J. Hedstrom and E. Houston, \emph{Pseudo-valuation domains II}, Houston J. Math. \textbf{4} (1978) 199--207.
\bibitem{HP}    W. Heinzer and I. Papick , \emph{The Radical Trace Property}, J. Algebra \textbf{112} (1988) 110--121.
\bibitem{HKLM1} E. Houston, S. Kabbaj, T. Lucas and A. Mimouni, \emph{Duals of Ideals in Pullback Constructions}, Lecture Notes Pure Appl. Math. Dekker \textbf{171} (1995) 263--276.
\bibitem{HuP}   J. A. Huckaba and I. J. Papick, \emph{When the dual of an ideal is a ring}, Manuscripta Math. \textbf{37} (1982) 67--85.
\bibitem{G1}    R. Gilmer, \emph{Some finiteness conditions on the set of overrings of an integral domain}, Proc. Amer. Math. Soc. \textbf{131} (2002) 2337--2346.
\bibitem{G2}    R. Gilmer, \emph{Multiplicative ideal theory}, Pure and Applied Mathematics, No. 12. Marcel Dekker, Inc., New York, 1972.
\bibitem{Ja1}   A. Jaballah, \emph{A lower bound for the number of intermediate rings}, Comm. Algebra \textbf{27} (1999) 1307--1311.
\bibitem{Ja2}   A. Jaballah, \emph{Finiteness of the set of intermediary rings in normal pairs}, Saitama Math. J. \textbf{17} (1999) 59--61.
\bibitem{Ja3}   A. Jaballah, \emph{The number of overrings of an integrally closed domain}, Expo. Math. \textbf{23} (2005) 353--360.
\bibitem{Ka}    I. Kaplansky, \emph{Commutative rings}, The University of Chicago Press, Chicago, 1974.
\bibitem{Kan}   B. G. Kang, \emph{${\star}$-Operations in integral domains}, PhD thesis, The University of Iowa, Iowa City, 1987.
\bibitem{M}     A. Mimouni, \emph{Integral domains in which each ideal is a $w$-ideal}, Comm. Algebra \textbf{33} (2005) 1345-1355.
\bibitem{M1}    A. Mimouni, \emph{Semistar-operations of finite character on integral domains}, J. Pure Appl. Algebra \textbf{200} (2005) 37-50.
\bibitem{MS}    A. Mimouni and M. Samman, \emph{semistar-operations on Valuation domains}, Internat. J. Comm. Rings \textbf{2(3)}, (2003) 131-141.
\bibitem{PT}    G. Picozza and F. Tartarone, \emph{When the semistar operation $\tilde{\star}$ is the identity}, Comm. Algebra to appear.
\bibitem{Z1}    M. Zafrullah, \emph {What $v$-coprimality can do for you, a survey article, to appear in ``Multiplicative Ideal Theory: a tribute to the work
                   of Robert Gilmer}," Eds:  J. W. Brewer, S. Glaz, W. J. Heinzer, and B. Olberding, Springer Science+Business Media, 2006.
\end{thebibliography}
\end{document}